\documentclass[reqno]{amsart}
\usepackage{amsmath}

\def\eqn#1{(\ref{eq:#1})}
\newcommand{\qBin}[3]{\genfrac{[}{]}{0pt}{0}{#1}{#2}_{#3}}

\newcommand{\qBinpwr}[4]{\genfrac{[}{]}{0pt}{0}{#1}{#2}_{#3}^{#4}}

\DeclareMathOperator{\mymod}{mod}
\theoremstyle{plain}
\newtheorem{thm}{Theorem}
\theoremstyle{remark}
\newtheorem*{rem}{Remark}

\numberwithin{equation}{section}
\allowdisplaybreaks[2]

\begin{document}

\title[Lattice paths, $q$-multinomials ...]
	{Lattice paths, $\boldsymbol{q}$-multinomials and
	 two variants of the Andrews-Gordon Identities}
\author[A.~Berkovich]{Alexander Berkovich}
\address{Department of Mathematics, The Pennsylvania State University,
         University Park, PA~16802, USA}
\email{alexb@math.psu.edu}
\thanks{This research was partially supported by SFB-grant
        F1305 of the Austrian FWF}
\author[P.~Paule]{Peter Paule}
\address{Research Institute for Symbolic Computation,
         Johannes Kepler University,
	 A--4040 Linz, Austria}
\email{Peter.Paule@risc.uni-linz.ac.at}

\begin{abstract}
A few years ago Foda, Quano, Kirillov and Warnaar proposed
and proved various finite analogs of the celebrated Andrews-Gordon identities.
In this paper we use these polynomial identities along with the combinatorial
techniques introduced in our recent paper to derive Garrett, Ismail, Stanton
type formulas for two variants of the Andrews-Gordon identities.
\end{abstract}

\maketitle

\section{Background and the first variant of the\\ Andrews-Gordon identities}
\label{sec:1}

\noindent
In 1961, Gordon~\cite{G} found a natural generalization of the Rogers-Ramanujan
partition theorem.
\begin{thm}\textup{(}Gordon\textup{)}
For all $\nu\ge 1$, $0\le s\le\nu$, the partitions of $N$ of the
frequency form $N=\sum_{j\ge 1}j f_j$ with $f_1\le s$ and $f_j+f_{j+1}\le\nu$,
$f_j\ge 0$ \textup{(}for all $j\ge 1$\textup{)} are equinumerous with the partitions of $N$ into
parts not congruent to $0$ or $\pm(s+1)$ modulo $2\nu+3$.
\label{thm:1}
\end{thm}

\noindent
Thirteen years later, Andrews~\cite{A1} proposed and proved the following analytic
counterpart to Gordon's theorem:
\begin{thm}\textup{(}Andrews\textup{)}
For all $\nu,s$ as in Theorem~\textup{\ref{thm:1}},  and $|q|<1$,
\begin{align}
\sum_{n_1,n_2,\ldots,n_\nu\ge 0} & \frac{q^{N_1^2+\ldots+N_\nu^2+N_{s+1}+\cdots+N_\nu}}
{(q)_{n_1}(q)_{n_2}\ldots(q)_{n_\nu}} = \frac{1}{(q)_\infty}\sum_{j=-\infty}^\infty
(-1)^j q^{\frac{j((2\nu+3)(j+1)-2(s+1))}{2}}\notag\\
& = \prod_{\stackrel{j\ge 1}{j\not\equiv 0,\pm (s+1)(\mymod 2\nu+3)}}\frac{1}{1-q^j},
\label{eq:1.1}
\end{align}
where
\begin{align}
N_i=\begin{cases}n_i+n_{i+1}+\cdots+n_\nu, & \mbox{if  } 1\le i\le \nu, \\
0, & \mbox{if  } i=\nu+1, \end{cases}
\label{eq:1.2}
\end{align}
and
\begin{equation}
(a;q)_\infty = (a)_\infty = \prod_{j\ge 0}(1-aq^j),
\label{eq:1.3}
\end{equation}
\begin{equation}
(a;q)_m = (a)_m = \prod_{j\ge 0}\frac{(1-aq^j)}{(1-aq^{j+m})}.
\label{eq:1.4}
\end{equation}
\label{thm:2}
\end{thm}

\noindent
We note that the last equality in \eqn{1.1} follows from Jacobi's triple
product identity~\cite[(II.28)]{GR}.

\begin{figure}[ht]
\begin{center}

\setlength{\unitlength}{0.5cm}

\begin{picture}(24,7.4)(-9,-1)
\small

\linethickness{0.2pt}
\put(-9,0){\vector(1,0){24}}
\put(0,-1){\vector(0,1){7}}

\put(-8,.15){\line(0,-1){.3}}
\put(-7,.15){\line(0,-1){.3}}
\put(-3,.15){\line(0,-1){.3}}
\put(-1,.15){\line(0,-1){.3}}
\put(1,.15){\line(0,-1){.3}}
\put(3,.15){\line(0,-1){.3}}
\put(7,.15){\line(0,-1){.3}}
\put(11,.15){\line(0,-1){.3}}
\put(13,.15){\line(0,-1){.3}}
\multiput(-.15,1)(0,1){4}{\line(1,0){.3}}

\put(-8,-.7){\makebox(0,0)[b]{$-8$}}
\put(-7,-.7){\makebox(0,0)[b]{$-7$}}
\put(-3,-.7){\makebox(0,0)[b]{$-3$}}
\put(-1,-.7){\makebox(0,0)[b]{$-1$}}
\put(1,-.7){\makebox(0,0)[b]{$1$}}
\put(3,-.7){\makebox(0,0)[b]{$3$}}
\put(7,-.7){\makebox(0,0)[b]{$7$}}
\put(11,-.7){\makebox(0,0)[b]{$11$}}
\put(13,-.7){\makebox(0,0)[b]{$13$}}
\put(15,-.7){\makebox(0,0)[b]{$x$}}

\put(0.3,1){\makebox(0,0)[l]{$1$}}
\put(0.3,2){\makebox(0,0)[l]{$2$}}
\put(0.3,3){\makebox(0,0)[l]{$3$}}
\put(0.3,4){\makebox(0,0)[l]{$4$}}
\put(0.3,6){\makebox(0,0)[l]{$y$}}

\multiput(-7.875,1)(.5,0){4}{\line(1,0){.25}}
\multiput(-1.875,1)(.5,0){4}{\line(1,0){.25}}
\multiput(2.125,1)(.5,0){4}{\line(1,0){.25}}

\thicklines
\put(-8,1){\line(1,1){1}}
\put(-7,2){\line(1,-1){2}}
\put(-5,0){\line(1,1){2}}
\put(-3,2){\line(1,-1){1}}
\put(-2,1){\line(1,1){1}}
\put(-1,2){\line(1,-1){2}}
\put(1,0){\line(1,1){2}}
\put(3,2){\line(1,-1){1}}
\put(4,1){\line(1,1){3}}
\put(7,4){\line(1,-1){4}}
\put(11,0){\line(1,0){2}}

\end{picture}

\caption{\label{fig:1}A Bressoud path starting at $(-8,1)$
and ending at $(13,0)$. Its weight is $-7-3-1+3+7=-1$.}

\end{center}
\end{figure}
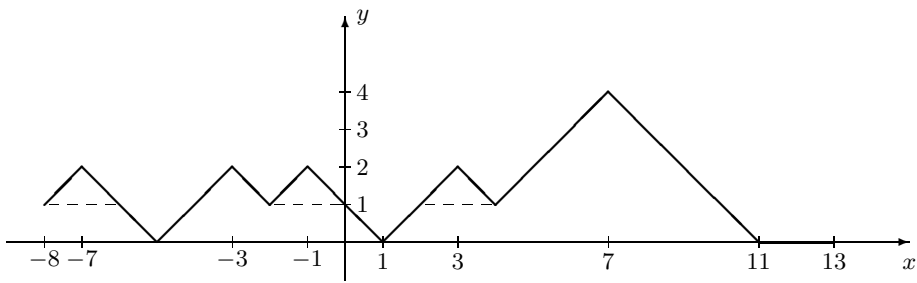

\noindent
Subsequently, Bressoud~\cite{B} interpreted the l.h.s.\ of \eqn{1.1} in terms of
weighted lattice paths. Bressoud path is made of three basic steps (see Fig.~\ref{fig:1}):
\begin{equation*}
\begin{array}{lll}
\mbox{NE step:} & \mbox{from } (i,j) \mbox{ to } (i+1,j+1), & \\
\mbox{SE step:} & \mbox{from } (i,j) \mbox{ to } (i+1,j-1), & \mbox{only allowed if } j>0, \\
\mbox{Horizontal step:} & \mbox{from } (i,0) \mbox{ to } (i+1,0), 
& \mbox{only allowed along } x \mbox{-axis},
\end{array}
\end{equation*}
with $i\in\mathbb{Z}$ and $j\in\mathbb{Z}_{\ge 0}$.

\noindent
To calculate the weight of Bressoud path we define the peak of a path as a
vertex preceded by the NE step and followed by the SE step. The height of a peak
is its $y$-coordinate, the weight of a peak is its $x$-coordinate. The weight $w(p)$
of path $p$ is defined as the sum of the weights of its peaks. In the above example
(Fig.~\ref{fig:1}), the path has five peaks:
$(-7,2),(-3,2),(-1,2),(3,2),(7,4)$ and its weight is $-7-3-1+3+7=-1$. The relative
height of a peak $(i,j)$ is the largest positive integer $h$, for which we can find two vertices
on the path: $(i',j-h),(i'',j-h)$ such that $i'<i<i''$ and such that between these two
vertices there are no peaks of height $>j$ and every peak of height $=j$ has weight $\ge i$.
The peaks in the above example (Fig.~\ref{fig:1}) have relative heights $1,2,1,1,4$, respectively.

\medskip
\noindent
We can now state Bressoud's result~\cite{B}.
\begin{thm} \textup{(}Bressoud\textup{)}
For $\nu,s$ as in Theorem~\textup{\ref{thm:1}} and $n_i\ge 0$, $1\le i\le\nu$
\begin{equation}
\frac{q^{N_1^2+\cdots+N_\nu^2+N_{s+1}+\cdots+N_\nu}}
{(q)_{n_1}(q)_{n_2}\ldots(q)_{n_\nu}} = \lim_{L\rightarrow\infty}
\sum_{p\in\mathbb P_{0,L}^\nu(\nu-s,0,\mathbf{n})} q^{w(p)},
\label{eq:1.5}
\end{equation}
where
\begin{equation}
\mathbf{n}=(n_1,n_2,\ldots,n_\nu),
\label{eq:1.6}
\end{equation}
and $\mathbb P_{x_1,x_2}^\nu(y_1,y_2,\mathbf n)$ denotes a collection of all Bressoud
lattice paths that start at $(x_1,y_1)$ and end at $(x_2,y_2)$, which have no peaks higher than
$\nu$, and, in addition, the number of peaks of relative height $j$ is $n_j$ for $1\le j\le\nu$.
\label{thm:3}
\end{thm}

\noindent
Actually, it is straightforward to refine Bressoud's analysis to show that
\begin{equation}
q^{N_1^2+\cdots+N_\nu^2+N_{s+1}+\cdots+N_\nu}
\prod_{i=1}^\nu \qBin{n_i+L-2\sum_{l=1}^i N_l-\alpha_{i,s}}{n_i}{q}
=\sum_{p\in\mathbb P_{0,L}^\nu(\nu-s,0,\mathbf{n})} q^{w(p)},
\label{eq:1.7}
\end{equation}
where
\begin{equation}
\alpha_{i,s}=\max(0,i-s),
\label{eq:1.8}
\end{equation}
and the $q$-binomial coefficients are defined as
\begin{align}
\qBin{n+m}{n}{q} =
\begin{cases}\frac{(q)_{n+m}}{(q)_n(q)_m}, & \mbox{for } n,m\in\mathbb{Z}_{\ge0}, \\
0, & \mbox{otherwise.} \end{cases}
\label{eq:1.9}
\end{align}
This leads to
\begin{align}
\sum_{\mathbf n}q^{N_1^2+\cdots+N_\nu^2+N_{s+1}+\cdots+N_\nu}
& \prod_{i=1}^L \qBin{n_i+L-2\sum_{l=1}^iN_l-\alpha_{i,s}}{n_i}{q}
& =C_{0,L}^\nu(\nu-s,0,q),
\label{eq:1.10}
\end{align}
where for $0\le s$, $b\le\nu$
\begin{equation}
C_{M,L}^\nu(s,b,q)=C_{M,L}^\nu(s,b):=
\sum_{p\in\tilde{\mathbb P}_{M,L}^\nu(s,b)} q^{w(p)}
\label{eq:1.11}
\end{equation}
and
\begin{equation}
\tilde{\mathbb P}_{M,L}^\nu(s,b)=
\sum_{\mathbf n}\mathbb P_{M,L}^\nu(s,b,\mathbf n).
\label{eq:1.12}
\end{equation}

\noindent
Here and in the following, summation over $\bf n$ is over all non-negative integer tuples
${\bf n}=(n_1,\ldots,n_{\nu})$.

\begin{rem}
We note that one can extract identity \eqn{1.7} from Lemma~4 in \cite{B3} with
$a=s+1$, $k=\nu$, $m_i=N_i$, $n_i=L-\alpha_{i,s}$, if one recognizes that, in this case, the
somewhat non-trivial conditions imposed on the peaks therein are equivalent to our concise
statement that all paths end at $(L,0)$.
\end{rem}

\noindent
Making use of \eqn{1.11}, one easily derives the following recursion relations
\begin{equation}
C_{0,L}^\nu(s,\nu)=C_{0,L-1}^\nu(s,\nu-1),
\label{eq:1.13}
\end{equation}
\begin{align}
C_{0,L}^\nu(s,b) & =C_{0,L-1}^\nu(s,b-1+\delta_{b,0})+C_{0,L-1}^\nu(s,b+1) \notag\\
& +(q^{L-1}-1)C_{0,L-2}^\nu(s,b), \quad 0\le b<\nu,
\label{eq:1.14}
\end{align}
and verifies the initial conditions
\begin{equation}
C_{0,0}^\nu(s,b)=\delta_{s,b},
\label{eq:1.15}
\end{equation}
where the Kronecker delta function $\delta_{i,j}$ is defined, as usual, as
\begin{align}
\delta_{i,j}=
\begin{cases} 1, & \mbox{if } i=j, \\
0, & \mbox{if } i\neq j. \end{cases}
\label{eq:1.16}
\end{align}
We note that formulas \eqn{1.13}--\eqn{1.15} specify the polynomials $C_{0,L}^\nu(s,b)$ uniquely.

\noindent
Next, for $L\equiv s+b(\mymod 2)$, we define polynomials $B_{s,b}^\nu(L,q)$ as
\begin{align}
B_{s,b}^\nu(L,q):=B_{s,b}^\nu(L) &:=\sum_{j=-\infty}^\infty \Bigg\{q^{j((2j+1)(2\nu+3)-2s)}
\qBin{L}{\frac{L+b-s}{2}+j(2\nu+3)}{q} \notag \\
& \qquad\qquad - q^{(2j+1)((2\nu+3)j+s)}\qBin{L}{\frac{L+b+s}{2}+j(2\nu+3)}{q}\Bigg\}.
\label{eq:1.17}
\end{align}
Employing the standard $q$-binomial recurrences \cite[(I.45)]{GR} one finds that
\begin{equation}
B_{s,b}^\nu(L)=B_{s,b-1}^\nu(L-1)+B_{s,b+1}^\nu(L-1)+(q^{L-1}-1)B_{s,b}^\nu(L-2).
\label{eq:1.18}
\end{equation}
It is not difficult to check that
\begin{equation}
B_{s,\nu+2}^\nu(L)=B_{2\nu+3-s,\nu+1}^\nu(L),
\label{eq:1.19}
\end{equation}
\begin{equation}
B_{s,1}^\nu(L)=B_{s,2}^\nu(L-1),
\label{eq:1.20}
\end{equation}
and for $1\le s,b\le\nu+1$
\begin{equation}
B_{s,b}^\nu(0)=\delta_{s,b}.
\label{eq:1.21}
\end{equation}
The formulas \eqn{1.13}--\eqn{1.15} and \eqn{1.18}--\eqn{1.21} imply that
\begin{align}
C_{0,L}^\nu(s,b)=
\begin{cases}B_{\nu+1-s,\nu+1-b}^\nu(L), & \mbox{if } L\equiv s+b(\mymod 2), \\
B_{\nu+2+s,\nu+1-b}^\nu(L), & \mbox{otherwise.} \end{cases}
\label{eq:1.22}
\end{align}
Indeed, both sides of \eqn{1.22} satisfy identical recurrences and initial conditions.

\noindent
Combining \eqn{1.10} and \eqn{1.22} we arrive at
\begin{thm}\textup{(}Foda, Quano and Kirillov\textup{)}
\begin{align}
\sum_{\mathbf n} q^{N_1^2+\cdots+N_\nu^2+N_s+\cdots+N_\nu}
& \prod_{i=1}^\nu \qBin{n_i+L-2\sum_{l=1}^i N_l-\alpha_{i,s-1}}{n_i}{q} \notag \\
& = \begin{cases}
B_{s,\nu+1}^\nu(L), & \mbox{if } L\not\equiv s+\nu(\mymod 2), \\
B_{2\nu+3-s,\nu+1}^\nu(L), & \mbox{otherwise.} \end{cases}
\label{eq:1.23}
\end{align}
\label{thm:4}
\end{thm}

\noindent
The above theorem was first proven in \cite{FQ} and \cite{K} in a somewhat different fashion.
Using the following limiting formulas
\begin{equation}
\lim_{L\rightarrow\infty}\qBin{L}{n}{q}=\frac{1}{(q)_n},
\label{eq:1.24}
\end{equation}
\begin{equation}
\lim_{L\rightarrow\infty}B_{s,b}^\nu(L)=\lim_{L\rightarrow\infty}B_{2\nu+3-s,b}^\nu(L)=
\prod_{\stackrel{j\ge 1}{j\not\equiv 0,\pm s(\mymod 2\nu+3)}} \frac{1}{1-q^j},
\label{eq:1.25}
\end{equation}
it is easy to check that in the limit $L\rightarrow\infty$, Theorem~\ref{thm:4} reduces to
Theorem~\ref{thm:2}.

\noindent
Recently, motivated by \cite{AKP,BM,GIS}, we investigated in \cite{BP}
the following multisums
\begin{equation}
\sum_{\mathbf{n}}\frac{q^{N_1^2+\cdots+N_\nu^2+N_i+\cdots+N_\nu-MN_1}}
{(q)_{n_1}(q)_{n_2}\ldots(q)_{n_\nu}}
\label{eq:1.26}
\end{equation}
with $1\le i\le\nu+1$, $M\in\mathbb Z_{\ge 0}$. The combinatorial analysis of \eqn{1.26}
with $\nu= 1$ given in \cite{BP} can be upgraded to the more general case $\nu\ge 1$ as
summarized in three steps below.

\medskip
\noindent
{\it First step}. We observe that
\begin{equation}
w(p')=w(p)-MN_1,
\label{eq:1.27}
\end{equation}
where the path $p'\in\mathbb P_{-M,L-M}^\nu(s,b,\mathbf n)$ is obtained from the path
$p\in\mathbb P_{0,L}^\nu(s,b,\mathbf n)$ by moving $p$ by $M$ units to the left along
the $x$-axis. Next, using \eqn{1.7} and \eqn{1.27} we derive that
\begin{align}
C_{-M,L-M}^\nu(0,0,q) & =
\sum_{\mathbf n}\sum_{p\in\mathbb P_{-M,L-M}^\nu(0,0,\mathbf n)}q^{w(p)}=
\sum_{\mathbf n}\sum_{p\in\mathbb P_{0,L}^\nu(0,0,\mathbf n)}q^{w(p)-MN_1} \notag \\
&=\sum_{\mathbf n}q^{N_1^2+\cdots+N_\nu^2-MN_1}
\prod_{i=1}^\nu \qBin{n_i+L-2\sum_{l=1}^i N_l}{n_i}{q}.
\label{eq:1.28}
\end{align}

\noindent
{\it Second step}. For $0\le M\le L$ every path $p\in\mathbb P_{-M,L-M}^\nu(0,0)$ consists of
two pieces joined together at some point $(0,s)$ with $0\le s\le\nu$. The first piece
belongs to $\mathbb P_{-M,0}^\nu(0,s)$ and the second one to $\mathbb P_{0,L-M}^\nu(s,0)$.
This observation is equivalent to
\begin{equation}
C_{-M,L-M}^\nu(0,0,q) =
\sum_{s=0}^\nu C_{-M,0}^\nu(0,s,q)\; C_{0,L-M}^\nu(s,0,q).
\label{eq:1.29}
\end{equation}
{\it Third step}.
\begin{equation}
C_{-M,0}^\nu(0,s,q)=
\sum_{p\in\mathbb P_{-M,0}^\nu(0,s)}q^{w(p)}=
\sum_{p\in\mathbb P_{0,M}^\nu(s,0)}q^{-w(p)}= C_{0,M}^\nu(s,0,\frac{1}{q}).
\label{eq:1.30}
\end{equation}
Combining \eqn{1.28}--\eqn{1.30} one obtains
\begin{align}
\sum_{\mathbf n}q^{N_1^2+\cdots+N_\nu^2-MN_1}
& \prod_{i=1}^\nu \qBin{n_i+L-2\sum_{l=1}^i N_l}{n_i}{q} \notag \\
& =  \sum_{s=0}^\nu C_{0,M}^\nu(s,0,\frac{1}{q}) \; C_{0,L-M}^\nu(s,0,q).
\label{eq:1.31}
\end{align}
Finally, letting $L$ tend to infinity, we find with the aid of \eqn{1.22},
\eqn{1.24}, \eqn{1.25} and \eqn{1.31} our first variant of the Andrews-Gordon
identities
\begin{align}
\sum_{\mathbf n}\frac{q^{N_1^2+\cdots+N_\nu^2-MN_1}}{(q)_{n_1}(q)_{n_2}\ldots(q)_{n_\nu}}
& = \sum_{\stackrel{s=1}{M\not\equiv\nu+s(\mymod 2)}}^{2\nu+2}
\frac{B_{s,\nu+1}^\nu(M,\frac{1}{q})}{\prod_{j\ge 1,j\not\equiv 0,\pm s
(\mymod 2\nu+3)}(1-q^j)}.
\label{eq:1.32}
\end{align}
Formula \eqn{1.32} is a special case of (3.21) in \cite{BP} with $s=\nu+1$. The other
cases there can be treated in a completely analogous manner.

\medskip
\noindent
Actually, neither the polynomial analogs \eqn{1.23} of \eqn{1.1}, nor the path interpretation
\eqn{1.10} of the Andrews-Gordon identities are unique. In particular, in \cite{W} Warnaar
considered the path space that is based on Gordon frequency conditions in Theorem~\ref{thm:1},
with an additional constraint that $f_j=0$ for $j>L$. This led him to the new polynomial
versions of the Andrews-Gordon identities. In the next section of this paper we will make 
essential use of Warnaar's analysis to investigate the following multisums
\begin{equation}
\sum_{\mathbf n}
\frac{q^{N_1^2+\cdots+N_\nu^2+N_s+\cdots+N_\nu-M(N_1+N_2+\cdots+N_\nu)}}
{(q)_{n_1}(q)_{n_2}\ldots(q)_{n_\nu}}, \quad M\in\mathbb Z_{\ge 0}.
\label{eq:1.33}
\end{equation}
The rest of this article is organized as follows. In Section~\ref{sec:2}, we briefly discuss
$q$-multinomial coefficients and Warnaar's terminating versions of the Andrews-Gordon identities.
In Section~\ref{sec:3}, we review a particle interpretation of Gordon's frequency conditions 
given in~\cite{W}. In Section~\ref{sec:4}, by following easy steps being similar to three 
steps above, we derive our main formulas for
\eqn{1.33} and their finite analogs. We conclude with a short description of prospects for
future work opened by this investigation.

\section{$q$-Multinomials and polynomial analogs of the\\ Andrews-Gordon identities}
\label{sec:2}

\noindent
We start by recalling the binomial theorem
\begin{equation}
(1+x)^L=\sum_{a=0}^L \binom{L}{a}x^a,
\label{eq:2.1}
\end{equation}
where $\binom{L}{a}$ is the usual binomial coefficient.

\noindent
By analogy, we introduce multinomial coefficients $\binom{L}{a}_\nu$ for
$a=0,1,\ldots,\nu L$ as the coefficients in the expansion 
\begin{equation}
(1+x+x^2+\cdots+x^\nu)^L=\sum_{a=0}^{\nu L} \binom{L}{a}_\nu x^a.
\label{eq:2.2}
\end{equation}
Multiple use of \eqn{2.1} yields an explicit sum representation
\begin{equation}
\binom{L}{a}_\nu =
\sum_{j_1+\cdots+j_\nu=a} \binom{L}{j_\nu}\binom{j_\nu}{j_{\nu-1}}\cdots\binom{j_2}{j_1}.
\label{eq:2.3}
\end{equation}
Building on the work of Andrews~\cite{A2}, Schilling~\cite{S} and Warnaar~\cite{W} have 
introduced the following $q$-analogs of \eqn{2.3}
\begin{equation}
\qBinpwr{L}{a}{\nu}{p}:=
\sum_{j_1+\cdots+j_\nu=a} q^{\sum_{l=2}^\nu j_{l-1}(L-j_l)-\sum_{l=1}^p j_l}
\qBin{L}{j_\nu}{q}\qBin{j_\nu}{j_{\nu-1}}{q}\cdots\qBin{j_2}{j_1}{q}
\label{eq:2.4}
\end{equation}
for $p=0,1,\ldots,\nu$.

\noindent
We list some important properties of $q$-multinomials \eqn{2.4}, which have been proven 
in \cite{S} and~\cite{W}.

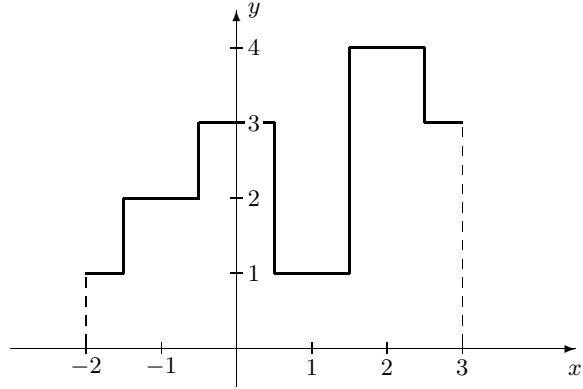
\begin{figure}[ht]
\begin{center}

\setlength{\unitlength}{0.5cm}

\begin{picture}(15,10.4)(-6,-1)
\small

\linethickness{0.2pt}
\put(-6,0){\vector(1,0){15}}
\put(0,-1){\vector(0,1){10}}

\multiput(-4,.15)(2,0){6}{\line(0,-1){.3}}
\multiput(-.15,2)(0,2){4}{\line(1,0){.3}}

\put(-4,-.7){\makebox(0,0)[b]{$-2$}}
\put(-2,-.7){\makebox(0,0)[b]{$-1$}}
\put(2,-.7){\makebox(0,0)[b]{$1$}}
\put(4,-.7){\makebox(0,0)[b]{$2$}}
\put(6,-.7){\makebox(0,0)[b]{$3$}}
\put(9,-.7){\makebox(0,0)[b]{$x$}}

\put(0.3,2){\makebox(0,0)[l]{$1$}}
\put(0.3,4){\makebox(0,0)[l]{$2$}}
\put(0.3,6){\makebox(0,0)[l]{$3$}}
\put(0.3,8){\makebox(0,0)[l]{$4$}}
\put(0.3,9){\makebox(0,0)[l]{$y$}}

\multiput(-4,0.125)(0,.5){4}{\line(0,1){.25}}
\multiput(6,0.125)(0,.5){12}{\line(0,1){.25}}

\thicklines
\put(-4,2){\line(1,0){1}}
\put(-3,2){\line(0,1){2}}
\put(-3,4){\line(1,0){2}}
\put(-1,4){\line(0,1){2}}
\put(-1,6){\line(1,0){1.2}}
\put(0.7,6){\line(1,0){0.3}}
\put(1,6){\line(0,-1){4}}
\put(1,2){\line(1,0){2}}
\put(3,2){\line(0,1){6}}
\put(3,8){\line(1,0){2}}
\put(5,8){\line(0,-1){2}}
\put(5,6){\line(1,0){1}}

\end{picture}

\caption{\label{fig:2}A lattice path representation of the sequence
$\{f_{-2}=1,f_{-1}=2,f_0=3,f_1=1,f_2=4,f_3=3\}$.
Here $M=-2$, $L=3$, $\nu=7$ and weight of the path
$=-1\cdot 2+0 \cdot 3+1 \cdot 1+2 \cdot 4=7$.}

\end{center}
\end{figure}

\medskip
\noindent
{\it {\bf Symmetries:}}
\begin{equation}
\qBinpwr{L}{a}{\nu}{p}=q^{(\nu-p)L-a}\qBinpwr{L}{\nu L-a}{\nu}{\nu-p}
\label{eq:2.5}
\end{equation}
and
\begin{equation}
\qBinpwr{L}{a}{\nu}{0} = \qBinpwr{L}{\nu L-a}{\nu}{0}.
\label{eq:2.6}
\end{equation}
{\it {\bf Recurrences:}}
\begin{equation}
\qBinpwr{L}{a}{\nu}{p} =
\sum_{m=0}^{\nu-p} q^{m(L-1)}\qBinpwr{L-1}{a-m}{\nu}{m}+
\sum_{m=\nu-p+1}^\nu q^{L(\nu-p)-m}\qBinpwr{L-1}{a-m}{\nu}{m}.
\label{eq:2.7}
\end{equation}
{\it {\bf $\boldsymbol{q}$-Deformed tautologies:}}
\begin{equation}
\qBinpwr{L}{a}{\nu}{p}+q^L\qBinpwr{L}{\nu L-a-p-1}{\nu}{p+1}=
q^L\qBinpwr{L}{a}{\nu}{p+1}+\qBinpwr{L}{\nu L-a-p-1}{\nu}{p}
\label{eq:2.8}
\end{equation}
where $p=-1,0,\ldots,\nu-1$ and
\begin{equation}
\qBinpwr{L}{a}{\nu}{-1}=0.
\label{eq:2.9}
\end{equation}
{\it {\bf Limiting behavior:}}
For $\frac{\nu L}{2}-A=0,1,\ldots,\nu L$
\begin{equation}
\lim_{L\rightarrow\infty}\qBinpwr{L}{\frac{\nu L}{2}-A}{\nu}{p}=
\begin{cases}
\frac{1}{(q)_\infty}, & \mbox{if } 0\le p<\frac{\nu}{2}, \\
\frac{1+q^A}{(q)_\infty}, & \mbox{if } \nu\equiv 0(\mymod 2) \mbox{ and } p=\frac{\nu}{2}, \\
\mbox{no limit}, & \mbox{if } \frac{\nu}{2}<p\le\nu.
\end{cases}
\label{eq:2.10}
\end{equation}
{\it {\bf Special values:}}
\begin{equation}
\qBinpwr{0}{a}{\nu}{p}=\delta_{a,0},
\label{eq:2.11}
\end{equation}
and
\begin{equation}
\qBinpwr{L}{a}{\nu}{p}=0, \quad \mbox{if } a<0 \mbox{ or } a>\nu L.
\label{eq:2.12}
\end{equation}
Next, for each $\nu\in\mathbb Z_{>0}$ we consider the ordered sequences of integers
$\{f_M,f_{M+1},\ldots,\linebreak[1]f_{L-1},f_L\}$ subject to Gordon's conditions
\begin{equation}
\begin{cases}
f_j\ge 0, & M\le j\le L, \\
f_j+f_{j+1}\le\nu, & M\le j< L.
\end{cases}
\label{eq:2.13}
\end{equation}
Each of these Gordon sequences can be represented graphically by a lattice path as illustrated 
by the example shown in Fig.~\ref{fig:2}.

\noindent
The weight $wg(p)$ of Gordon path $p$ described above is
\begin{equation}
wg(p)=\sum_{j=M+1}^{L-1}jf_j.
\label{eq:2.14}
\end{equation}
For $0\le s$, $b\le\nu$, $L,M\in\mathbb Z$ we perform a weighted path count with the help
of the following polynomials
\begin{equation}
G_{M,L}^\nu(s,b,q):=\sum_{p\in\tilde{\mathbb S}_{M,L}^\nu(s,b)}q^{wg(p)},
\label{eq:2.15}
\end{equation}
where $\tilde{\mathbb S}_{M,L}^\nu(s,b)$ denotes the space of all Gordon paths
subject to \eqn{2.13} that start at $(M,s)$ and end at $(L,b)$.

\noindent
In \cite{W}, Warnaar, building on the work of Andrews and Baxter~\cite{AB}, showed
in a recursive fashion that for $0\le s$, $b\le\nu$
\begin{equation}
G_{0,L}^\nu(\nu-s,b,q)=
\begin{cases}
W_{s+1,b}^\nu(L,q), & \mbox{if } b+s\equiv\nu(L+1)(\mymod 2), \\
W_{2\nu+2-s,b}^\nu(L,q), & \mbox{if } b+s\not\equiv\nu(L+1)(\mymod 2),
\end{cases}
\label{eq:2.16}
\end{equation}
where
\begin{align}
W_{s,b}^\nu(L,q) & := \sum_{j=-\infty}^\infty \Bigg\{
q^{j((2j+1)(2\nu+3)-2s)}
\qBinpwr{L}{\frac{\nu(L+1)-s-b+1}{2}+(2\nu+3)j}{\nu}{b} \notag \\
& \qquad\qquad - q^{(2j+1)((2\nu+3)j+s)}
\qBinpwr{L}{\frac{\nu(L+1)+s-b+1}{2}+(2\nu+3)j}{\nu}{b} \Bigg\}.
\label{eq:2.17}
\end{align}
More specifically, it was proven in \cite{W} that both sides of \eqn{2.16}
satisfy the same recurrences
\begin{equation}
G_{0,L}^\nu(\nu-s,b,q)=\sum_{l=0}^{\nu-b}q^{(L-1)l}G_{0,L-1}^\nu(\nu-s,l,q)
\label{eq:2.18}
\end{equation}
and the same initial conditions
\begin{equation}
G_{0,0}^\nu(\nu-s,b,q)=\delta_{s+b,\nu}.
\label{eq:2.19}
\end{equation}
Also in \cite{W}, a particle interpretation of Gordon paths was given and, as
a result, another representation for $G_{0,L}^\nu(s,b)$ was obtained; namely,
\begin{equation}
G_{0,L}^\nu(s,b)=F_{\nu-s,\nu-b}^\nu(L,q),
\label{eq:2.20}
\end{equation}
where $0\le s$, $b\le\nu$, $L\ge 2$ and
\begin{equation}
F_{s,b}^\nu(L,q):=
\sum_{\mathbf n}q^{N_1^2+\cdots+N_\nu^2+N_{s+1}+\cdots+N_\nu}
\prod_{i=1}^\nu
\qBin{n_i+iL-2\sum_{l=1}^iN_l-\alpha_{i,s}-\alpha_{i,b}}{n_i}{q}.
\label{eq:2.21}
\end{equation}
Comparing \eqn{2.16} and \eqn{2.20} we arrive at the polynomial identities
\begin{equation}
F_{s,b}^\nu(L,q)=
\begin{cases}
W_{s+1,\nu-b}^\nu(L,q), & \mbox{if } b+s\equiv\nu L(\mymod 2), \\
W_{2\nu+2-s,\nu-b}^\nu(L,q), & \mbox{otherwise}, \end{cases}
\label{eq:2.22}
\end{equation}
which in the limit $L\rightarrow\infty$ reduce to the Andrews-Gordon identities \eqn{1.1}.
It is important to realize that while \eqn{1.23} and \eqn{2.22} are identical in the
limit $L\rightarrow\infty$, these identities are substantially different for finite $L$.

\begin{rem}
It should be noted that \eqn{2.20} with $s=b=0$ is a corollary of Bressoud's Lemma 3
in~\cite{B3}.  It appears that this Lemma, as stated, is true for $a=k+1$ only. For
$a=1,2,\dots,k$, the generating function
$c_{a,k}(n,2n,\ldots,kn;j)$ therein should be corrected to $c_{a,k}(n-\alpha_{1,a-1},
2n-\alpha_{2,a-1},\ldots,kn-\alpha_{k,a-1};j)$.
\end{rem}

\section{Gordon paths and lattice particles}
\label{sec:3}

\begin{figure}[ht]
\begin{center}

\setlength{\unitlength}{0.5cm}

\begin{picture}(24,11.4)(-2,-1)
\small

\linethickness{0.2pt}
\put(-2,0){\vector(1,0){24}}
\put(0,-1){\vector(0,1){11}}

\multiput(1,.15)(2,0){2}{\line(0,-1){.3}}
\put(20,.15){\line(0,-1){.3}}
\put(-.15,1){\line(1,0){.3}}
\put(-.15,7){\line(1,0){.3}}
\put(-.15,8){\line(1,0){.3}}

\put(-.5,-.3){\makebox(0,0)[t]{$M$}}
\put(1,-.3){\makebox(0,0)[t]{$M+1$}}
\put(3,-.3){\makebox(0,0)[t]{$M+3$}}
\put(20,-.3){\makebox(0,0)[t]{$L$}}
\put(22,-.4){\makebox(0,0)[t]{$x$}}

\put(-0.3,1){\makebox(0,0)[r]{$1$}}
\put(-0.3,7){\makebox(0,0)[r]{$\nu-1$}}
\put(-0.3,8){\makebox(0,0)[r]{$\nu$}}
\put(-0.3,10){\makebox(0,0)[r]{$y$}}

\thicklines
\put(0,0){\line(1,0){0.5}}
\put(0.5,0){\line(0,1){8}}
\put(0.5,8){\line(1,0){1}}
\put(1.5,8){\line(0,-1){8}}
\put(1.5,0){\line(1,0){1}}
\put(2.5,0){\line(0,1){8}}
\put(2.5,8){\line(1,0){1}}
\put(3.5,8){\line(0,-1){8}}
\put(3.5,0){\line(1,0){0.5}}
\put(5,0){\line(1,0){0.5}}
\put(5.5,0){\line(0,1){8}}
\put(5.5,8){\line(1,0){1}}
\put(6.5,8){\line(0,-1){8}}
\put(6.5,0){\line(1,0){1}}

\put(7.5,0){\line(0,1){7}}
\put(7.5,7){\line(1,0){1}}
\put(8.5,7){\line(0,-1){7}}
\put(8.5,0){\line(1,0){0.5}}
\put(10,0){\line(1,0){0.5}}
\put(10.5,0){\line(0,1){7}}
\put(10.5,7){\line(1,0){1}}
\put(11.5,7){\line(0,-1){7}}
\put(11.5,0){\line(1,0){0.5}}

\put(14.5,0){\line(1,0){0.5}}
\put(15,0){\line(0,1){1}}
\put(15,1){\line(1,0){1}}
\put(16,1){\line(0,-1){1}}
\put(16,0){\line(1,0){0.5}}
\put(17.5,0){\line(1,0){0.5}}
\put(18,0){\line(0,1){1}}
\put(18,1){\line(1,0){1}}
\put(19,1){\line(0,-1){1}}
\put(19,0){\line(1,0){1}}

\put(4.5,6){\makebox(0,0){$\dots$}}
\put(9.5,5){\makebox(0,0){$\dots$}}
\put(13.25,3){\makebox(0,0){$\dots$}}
\put(17,0.5){\makebox(0,0){$\dots$}}

\put(0.5,8.2){\makebox(6,1){$\overbrace{\hspace{3cm}}^{\textstyle n_\nu}$}}
\put(7.5,7.2){\makebox(4,1){$\overbrace{\hspace{2cm}}^{\textstyle n_{\nu-1}}$}}
\put(15,1.2){\makebox(4,1){$\overbrace{\hspace{2cm}}^{\textstyle n_1}$}}

\end{picture}

\caption{\label{fig:3}The minimal path in $\tilde{\mathbb{S}}_{M,L}^\nu(0,0)$
of particle content $\mathbf{n}=(n_1,n_2,\dots,n_\nu)$.}

\end{center}
\end{figure}
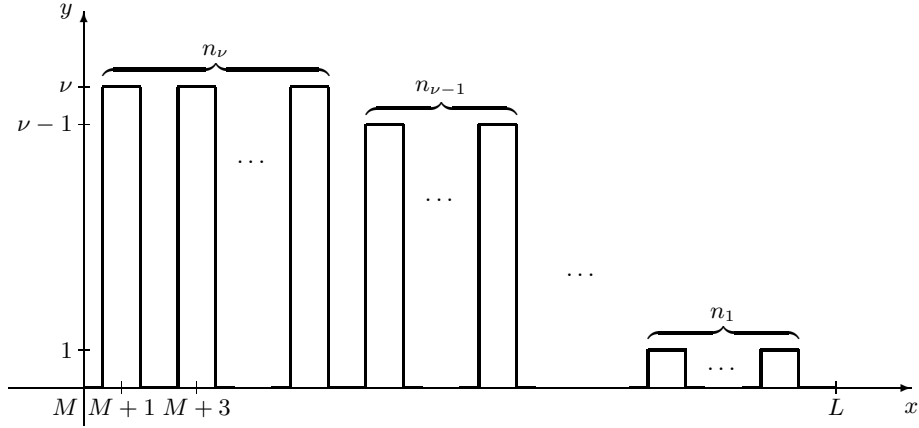

\noindent
We begin our particle description of Gordon lattice paths by considering first paths
in $\tilde{\mathbb S}_{M,L}^\nu(0,0)$. Following \cite{W}, we introduce a special kind of paths
from which all other paths in $\tilde{\mathbb S}_{M,L}^\nu(0,0)$ can be constructed. These
paths, termed the minimal paths in \cite{W}, are shown in Fig.~\ref{fig:3}.

\begin{figure}[ht]
\begin{center}

\setlength{\unitlength}{0.5cm}

\begin{picture}(21,5.9)(0,-0.8)
\small

\linethickness{0.2pt}
\put(0,0){\line(1,0){3.5}}
\put(5,0){\line(1,0){3.5}}
\put(10,0){\line(1,0){3.5}}
\put(17,0){\line(1,0){3.5}}

\multiput(1,.15)(2,0){2}{\line(0,-1){.3}}
\multiput(6,.15)(2,0){2}{\line(0,-1){.3}}
\multiput(11,.15)(2,0){2}{\line(0,-1){.3}}
\multiput(18,.15)(2,0){2}{\line(0,-1){.3}}

\put(1,-.8){\makebox(0,0)[b]{$j$}}
\put(3,-.8){\makebox(0,0)[b]{$j+2$}}
\put(6,-.8){\makebox(0,0)[b]{$j$}}
\put(8,-.8){\makebox(0,0)[b]{$j+2$}}
\put(11,-.8){\makebox(0,0)[b]{$j$}}
\put(13,-.8){\makebox(0,0)[b]{$j+2$}}
\put(18,-.8){\makebox(0,0)[b]{$j$}}
\put(20,-.8){\makebox(0,0)[b]{$j+2$}}

\thicklines
\put(0,0){\line(1,0){0.5}}
\put(0.5,0){\line(0,1){4.5}}
\put(0.5,4.5){\line(1,0){1}}
\put(1.5,4.5){\line(0,-1){4.5}}
\put(1.5,0){\line(1,0){0.5}}

\put(5,0){\line(1,0){0.5}}
\put(5.5,0){\line(0,1){4}}
\put(5.5,4){\line(1,0){1}}
\put(6.5,4){\line(0,-1){3.5}}
\put(6.5,0.5){\line(1,0){1}}
\put(7.5,0.5){\line(0,-1){0.5}}
\put(7.5,0){\line(1,0){0.5}}

\put(10,0){\line(1,0){0.5}}
\put(10.5,0){\line(0,1){3.5}}
\put(10.5,3.5){\line(1,0){1}}
\put(11.5,3.5){\line(0,-1){2.5}}
\put(11.5,1){\line(1,0){1}}
\put(12.5,1){\line(0,-1){1}}
\put(12.5,0){\line(1,0){0.5}}

\put(18,0){\line(1,0){0.5}}
\put(18.5,0){\line(0,1){4.5}}
\put(18.5,4.5){\line(1,0){1}}
\put(19.5,4.5){\line(0,-1){4.5}}
\put(19.5,0){\line(1,0){0.5}}

\put(1,4.8){\makebox(0,0)[b]{$t$}}
\put(6,4.3){\makebox(0,0)[b]{$t-1$}}
\put(7,0.8){\makebox(0,0)[b]{$1$}}
\put(11,3.8){\makebox(0,0)[b]{$t-2$}}
\put(12,1.3){\makebox(0,0)[b]{$2$}}
\put(19,4.8){\makebox(0,0)[b]{$t$}}

\put(3.5,2.5){\makebox(0,0)[b]{$\longrightarrow$}}
\put(8.5,2.5){\makebox(0,0)[b]{$\longrightarrow$}}
\put(15,2.5){\makebox(0,0)[b]{$\longrightarrow\dots\longrightarrow$}}

\end{picture}

\caption{\label{fig:4}A particle of charge $t$ in free motion from
$j$ to $j+1$.}

\end{center}
\end{figure}
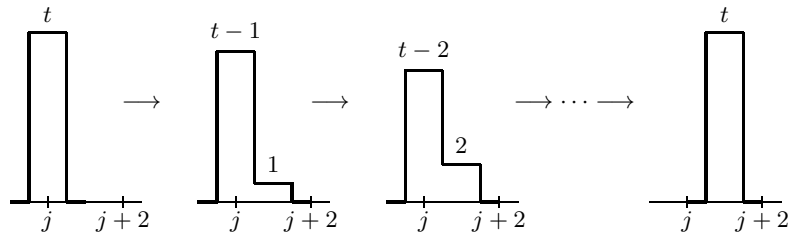

\noindent
In a minimal path, each column with a non-zero height $t$ $(0<t\le\nu)$ is interpreted
as a particle of charge $t$. Two adjacent particles in a minimal path are separated by
a single empty column. In order to construct an arbitrary non-minimal path in
$\mathbb{\tilde S}_{M,L}^\nu(0,0)$ out of one and only one minimal path, we need to introduce
the rules of particle motion from left to right. Fig.~\ref{fig:4} illustrates these rules in the
simplest case of an isolated particle of charge $t$ going from $j$ to $j+1$. The
complete transition requires $t$ elementary moves.

\noindent
Next, we consider the motion of a particle of charge $t$ through a path configuration
shown in Fig.~\ref{fig:5} and Fig.~\ref{fig:6}(a).

\begin{figure}[ht]
\begin{center}

\setlength{\unitlength}{0.5cm}

\begin{picture}(15,6.5)(0,-0.8)
\small

\linethickness{0.2pt}
\put(0,0){\line(1,0){15}}

\multiput(5,.15)(2,0){2}{\line(0,-1){.3}}

\put(5,-.8){\makebox(0,0)[b]{$j$}}
\put(7,-.8){\makebox(0,0)[b]{$j+2$}}

\thicklines
\put(4,0){\line(1,0){0.5}}
\put(4.5,0){\line(0,1){3.5}}
\put(4.5,3.5){\line(1,0){1}}
\put(5.5,3.5){\line(0,-1){3.5}}
\put(5.5,0){\line(1,0){1}}
\put(6.5,0){\line(0,1){4}}
\put(6.5,4){\line(1,0){1}}
\put(7.5,4){\line(0,1){1}}
\put(7.5,5){\line(1,0){1}}

\put(5,3.8){\makebox(0,0)[b]{$t$}}
\put(7,4.3){\makebox(0,0)[b]{$s$}}
\put(8,5.3){\makebox(0,0)[b]{$u$}}
\put(9.1,5){\makebox(0,0){$\dots$}}

\end{picture}

\caption{\label{fig:5}$t \le s$}

\end{center}
\end{figure}

\begin{figure}[ht]
\begin{center}

\setlength{\unitlength}{0.5cm}

\begin{picture}(21,6.7)(-1,-0.8)
\small

\linethickness{0.2pt}
\put(0,0){\line(1,0){6}}
\put(14,0){\line(1,0){6}}

\multiput(1,.15)(2,0){2}{\line(0,-1){.3}}
\multiput(15,.15)(2,0){2}{\line(0,-1){.3}}

\put(1,-.8){\makebox(0,0)[b]{$j$}}
\put(3,-.8){\makebox(0,0)[b]{$j+2$}}
\put(15,-.8){\makebox(0,0)[b]{$j$}}
\put(17,-.8){\makebox(0,0)[b]{$j+2$}}

\thicklines
\put(0,0){\line(1,0){0.5}}
\put(0.5,0){\line(0,1){5}}
\put(0.5,5){\line(1,0){1}}
\put(1.5,5){\line(0,-1){5}}
\put(1.5,0){\line(1,0){1}}
\put(2.5,0){\line(0,1){2.5}}
\put(2.5,2.5){\line(1,0){1}}

\put(14,0){\line(1,0){0.5}}
\put(14.5,0){\line(0,1){2.5}}
\put(14.5,2.5){\line(1,0){1}}
\put(15.5,2.5){\line(0,1){1.5}}
\put(15.5,4){\line(1,0){1}}
\put(16.5,4){\line(0,-1){1.5}}
\put(16.5,2.5){\line(1,0){1}}

\put(1,5.3){\makebox(0,0)[b]{$t$}}
\put(3,2.8){\makebox(0,0)[b]{$s$}}
\put(4.1,2.5){\makebox(0,0){$\dots$}}
\put(15,2.8){\makebox(0,0)[b]{$s$}}
\put(16,4.3){\makebox(0,0)[b]{$t-s$}}
\put(17,2.8){\makebox(0,0)[b]{$s$}}
\put(18.1,2.5){\makebox(0,0){$\dots$}}

\put(-1,5.5){\makebox(0,0)[l]{(a)}}
\put(13,5.5){\makebox(0,0)[l]{(b)}}

\put(7,3.5){\shortstack{After $(t-s)$\\ elementary moves}}
\put(8,2){\vector(1,0){3}}
\end{picture}

\caption{\label{fig:6}$t > s$}

\end{center}
\end{figure}

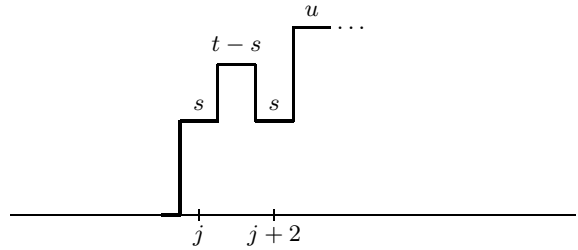
\begin{figure}[ht]
\begin{center}

\setlength{\unitlength}{0.5cm}

\begin{picture}(15,6.5)(0,-0.8)
\small

\linethickness{0.2pt}
\put(0,0){\line(1,0){15}}

\multiput(5,.15)(2,0){2}{\line(0,-1){.3}}

\put(5,-.8){\makebox(0,0)[b]{$j$}}
\put(7,-.8){\makebox(0,0)[b]{$j+2$}}

\thicklines
\put(4,0){\line(1,0){0.5}}
\put(4.5,0){\line(0,1){2.5}}
\put(4.5,2.5){\line(1,0){1}}
\put(5.5,2.5){\line(0,1){1.5}}
\put(5.5,4){\line(1,0){1}}
\put(6.5,4){\line(0,-1){1.5}}
\put(6.5,2.5){\line(1,0){1}}
\put(7.5,2.5){\line(0,1){2.5}}
\put(7.5,5){\line(1,0){1}}

\put(5,2.8){\makebox(0,0)[b]{$s$}}
\put(6,4.3){\makebox(0,0)[b]{$t-s$}}
\put(7,2.8){\makebox(0,0)[b]{$s$}}
\put(8,5.3){\makebox(0,0)[b]{$u$}}
\put(9.1,5){\makebox(0,0){$\dots$}}

\end{picture}

\caption{\label{fig:7}$u \ge t-s$}

\end{center}
\end{figure}

\noindent
In case of the path configuration in Fig.~\ref{fig:5}, the particle of charge $t$ can not
make any further move to the right. In the case shown in Fig.~\ref{fig:6}(a), we can make
$t-s$ elementary moves to obtain the new path configuration in Fig.~\ref{fig:6}(b). What
happens next depends on the height of the column at $j+3$, as illustrated
in Fig.~\ref{fig:7} and Fig.~\ref{fig:8}.

\noindent
In case of the path configuration in Fig.~\ref{fig:7}, the particle of charge $t$ can not
move any further.

\begin{figure}[ht]
\begin{center}

\setlength{\unitlength}{0.5cm}

\begin{picture}(21,5.8)(-1,-0.8)
\small

\linethickness{0.2pt}
\put(0,0){\line(1,0){6}}
\put(14,0){\line(1,0){6}}

\multiput(1,.15)(2,0){2}{\line(0,-1){.3}}
\multiput(15,.15)(2,0){2}{\line(0,-1){.3}}

\put(1,-.8){\makebox(0,0)[b]{$j$}}
\put(3,-.8){\makebox(0,0)[b]{$j+2$}}
\put(15,-.8){\makebox(0,0)[b]{$j$}}
\put(17,-.8){\makebox(0,0)[b]{$j+2$}}

\thicklines
\put(0,0){\line(1,0){0.5}}
\put(0.5,0){\line(0,1){2.5}}
\put(0.5,2.5){\line(1,0){1}}
\put(1.5,2.5){\line(0,1){1.5}}
\put(1.5,4){\line(1,0){1}}
\put(2.5,4){\line(0,-1){1.5}}
\put(2.5,2.5){\line(1,0){1}}
\put(3.5,2.5){\line(0,1){1}}
\put(3.5,3.5){\line(1,0){1}}

\put(14,0){\line(1,0){0.5}}
\put(14.5,0){\line(0,1){2.5}}
\put(14.5,2.5){\line(1,0){1}}
\put(15.5,2.5){\line(0,1){1}}
\put(15.5,3.5){\line(1,0){1}}
\put(16.5,3.5){\line(0,1){0.5}}
\put(16.5,4){\line(1,0){1}}
\put(17.5,4){\line(0,-1){0.5}}
\put(17.5,3.5){\line(1,0){1}}

\put(1,2.8){\makebox(0,0)[b]{$s$}}
\put(2,4.3){\makebox(0,0)[b]{$t-s$}}
\put(3,2.8){\makebox(0,0)[b]{$s$}}
\put(4,3.8){\makebox(0,0)[b]{$u$}}
\put(5.1,3.5){\makebox(0,0){$\dots$}}
\put(15,2.8){\makebox(0,0)[b]{$s$}}
\put(16,3.8){\makebox(0,0)[b]{$u$}}
\put(17,4.3){\makebox(0,0)[b]{$t-u$}}
\put(18,3.8){\makebox(0,0)[b]{$u$}}
\put(19.1,3.5){\makebox(0,0){$\dots$}}

\put(-1,4.5){\makebox(0,0)[l]{(a)}}
\put(13,4.5){\makebox(0,0)[l]{(b)}}

\put(7,2.5){\shortstack{After $(t-s-u)$\\ elementary moves}}
\put(8,1){\vector(1,0){3}}
\end{picture}

\caption{\label{fig:8}$u < t - s$}

\end{center}
\end{figure}

\noindent
In the case shown in Fig.~\ref{fig:8}(a) we can make $t-u-s$ elementary moves to end up with
the path configuration in Fig.~\ref{fig:8}(b). Ignoring the first column at $j$, we see that
the last configuration is practically the same as the one in Fig.~\ref{fig:6}(b) with $s$
replaced by $u$. It means that we can keep on moving the
particle of charge $t$ according to the rules in Fig.~\ref{fig:4}--Fig.~\ref{fig:8}.

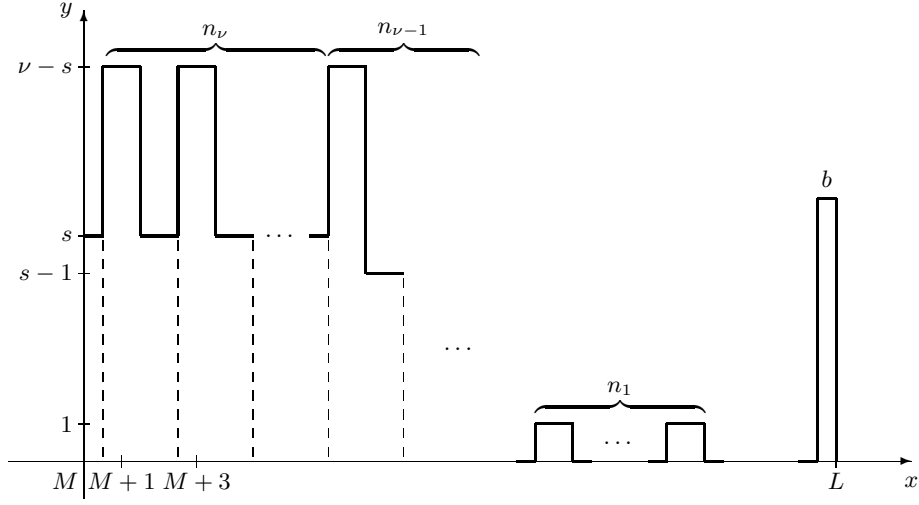
\begin{figure}[ht]
\begin{center}

\setlength{\unitlength}{0.5cm}

\begin{picture}(24,13.4)(-2,-1)
\small

\linethickness{0.2pt}
\put(-2,0){\vector(1,0){24}}
\put(0,-1){\vector(0,1){13}}

\multiput(1,.15)(2,0){2}{\line(0,-1){.3}}
\put(20,.15){\line(0,-1){.3}}
\put(-.15,1){\line(1,0){.3}}
\put(-.15,5){\line(1,0){.3}}
\put(-.15,6){\line(1,0){.3}}
\put(-.15,10.5){\line(1,0){.3}}

\put(-.5,-.3){\makebox(0,0)[t]{$M$}}
\put(1,-.3){\makebox(0,0)[t]{$M+1$}}
\put(3,-.3){\makebox(0,0)[t]{$M+3$}}
\put(20,-.3){\makebox(0,0)[t]{$L$}}
\put(22,-.4){\makebox(0,0)[t]{$x$}}

\put(-0.3,1){\makebox(0,0)[r]{$1$}}
\put(-0.3,5){\makebox(0,0)[r]{$s-1$}}
\put(-0.3,6){\makebox(0,0)[r]{$s$}}
\put(-0.3,10.5){\makebox(0,0)[r]{$\nu-s$}}
\put(-0.3,12){\makebox(0,0)[r]{$y$}}

\multiput(0.5,0.125)(0,0.5){12}{\line(0,1){.25}}
\multiput(2.5,0.125)(0,0.5){12}{\line(0,1){.25}}
\multiput(4.5,0.125)(0,0.5){12}{\line(0,1){.25}}
\multiput(6.5,0.125)(0,0.5){12}{\line(0,1){.25}}
\multiput(8.5,0.125)(0,0.5){10}{\line(0,1){.25}}

\thicklines
\put(0,6){\line(1,0){0.5}}
\put(0.5,6){\line(0,1){4.5}}
\put(0.5,10.5){\line(1,0){1}}
\put(1.5,10.5){\line(0,-1){4.5}}
\put(1.5,6){\line(1,0){1}}
\put(2.5,6){\line(0,1){4.5}}
\put(2.5,10.5){\line(1,0){1}}
\put(3.5,10.5){\line(0,-1){4.5}}
\put(3.5,6){\line(1,0){1}}
\put(6,6){\line(1,0){0.5}}
\put(6.5,6){\line(0,1){4.5}}
\put(6.5,10.5){\line(1,0){1}}
\put(7.5,10.5){\line(0,-1){5.5}}
\put(7.5,5){\line(1,0){1}}

\put(11.5,0){\line(1,0){0.5}}
\put(12,0){\line(0,1){1}}
\put(12,1){\line(1,0){1}}
\put(13,1){\line(0,-1){1}}
\put(13,0){\line(1,0){0.5}}
\put(15,0){\line(1,0){0.5}}
\put(15.5,0){\line(0,1){1}}
\put(15.5,1){\line(1,0){1}}
\put(16.5,1){\line(0,-1){1}}
\put(16.5,0){\line(1,0){0.5}}

\put(19,0){\line(1,0){0.5}}
\put(19.5,0){\line(0,1){7}}
\put(19.5,7){\line(1,0){0.5}}
\put(20,7){\line(0,-1){7}}

\put(5.25,6){\makebox(0,0){$\dots$}}
\put(10,3){\makebox(0,0){$\dots$}}
\put(14.25,0.5){\makebox(0,0){$\dots$}}

\put(19.75,7.3){\makebox(0,0)[b]{$b$}}

\put(0.5,10.7){\makebox(6,1)[b]{$\overbrace{\hspace{2.9cm}}^{\textstyle n_\nu}$}}
\put(6.5,10.7){\makebox(4,1)[b]{$\overbrace{\hspace{2cm}}^{\textstyle n_{\nu-1}}$}}
\put(12,1.2){\makebox(4.5,1){$\overbrace{\hspace{2.25cm}}^{\textstyle n_1}$}}

\end{picture}

\caption{\label{fig:9}The minimal path in $\tilde{\mathbb{S}}_{M,L}^\nu(s,b)$
of particle content $\mathbf{n}=(n_1,n_2,\dots,n_\nu)$. The dashed lines separate 
different particles. All particles are arranged according to their charges in 
non-decreasing order from left to right.}

\end{center}
\end{figure}

\noindent
Actually, it is easy to modify the above discussion in order to deal with the general
boundary conditions $0\le s$, $b\le\nu$. All we need to do is to refine the notion of a
minimal path as in Fig.~\ref{fig:9}.

\noindent
Two adjacent particles in Fig.~\ref{fig:9} are separated by at most one empty column. Two
half-columns at $M$ and $L$ are not interpreted as particles.

\noindent
In \cite{W}, Warnaar proved that using rules of particle motion, one can obtain
each non-minimal path from one and only one minimal path in a completely bijective
fashion. The particle content $\mathbf n$ of $p\in\tilde{\mathbb S}_{M,L}^\nu(s,b)$
can be determined by reducing $p$ to its minimal image $p_{\min}$.

\noindent
Now, since the sum of heights $\sum_{j=M+1}^{L-1}f_j$ of the path
$p\in\mathbb S_{M,L}^\nu(s,b,\mathbf n)$ is the invariant of the motion we find by
counting the heights of $p_{\min}$ that
\begin{equation}
\sum_{j=M-1}^{L-1}f_j=\sum_{j=1}^\nu jn_j=N_1+N_2+\cdots+N_\nu,
\label{eq:3.1}
\end{equation}
where $\mathbb S_{M,L}^\nu(s,b,\mathbf n)$ denotes the space of all Gordon paths of
the particle content $\mathbf n$, which start at $(M,s)$ and end at $(L,b)$.

\noindent
Equation \eqn{3.1} implies that
\begin{equation}
wg(\tilde p)=wg(p)-M(N_1+N_2+\cdots+N_\nu),
\label{eq:3.2}
\end{equation}
where $\tilde p\in\mathbb S_{-M,L-M}^\nu(s,b,\mathbf n)$ is obtained from
$p\in\mathbb S_{0,L}^\nu(s,b,\mathbf n)$ by moving $p$ by $M$ units to the left along
the $x$-axis.

\noindent
Formula \eqn{3.2} will play an important role in the sequel. We shall also require
another result established in~\cite{W}:
\begin{thm}\textup{(}Warnaar\textup{)}
\begin{align}
q^{N_1^2+\cdots+N_\nu^2+N_{s+1}+\cdots+N_\nu}
& \prod_{i=1}^\nu \qBin{n_i+iL-2\sum_{l=1}^i N_l-\alpha_{i,s}-\alpha_{i,b}}{n_i}{q} \notag\\
& = \sum_{p\in\mathbb S_{0,L}^\nu(\nu-s,\nu-b,\mathbf n)}q^{wg(p)},
\label{eq:3.3}
\end{align}
where $0\le s$, $b\le\nu$ and $L\ge 2$.
\label{thm:5}
\end{thm}

\noindent
Note that \eqn{2.20} is an immediate consequence of Theorem~\ref{thm:5}.

\section{The second variant of the Andrews-Gordon identities}
\label{sec:4}

\noindent
Having collected the necessary background information, we can derive the second variant of the
Andrews-Gordon identities by following three easy steps very similar to those taken to derive
our first variant in Section~\ref{sec:1}.

\medskip
\noindent
{\it First step}. We generalize \eqn{3.3} as
\begin{align}
\sum_{p\in\mathbb S_{-M,L-M}^\nu(\nu-s,\nu-b,\mathbf n)}q^{wg(p)}
& \stackrel{\mbox{\scriptsize{by \eqn{3.2}}}}{=}
\sum_{p\in\mathbb S_{0,L}^\nu(\nu-s,\nu-b,\mathbf n)}q^{wg(p)-M(N_1+\cdots+N_\nu)}\notag\\
& \stackrel{\mbox{\scriptsize{by \eqn{3.3}}}}{=}
q^{N_1^2+\cdots+N_\nu^2+N_{s+1}+\cdots+N_\nu-M(N_1+\cdots+N_\nu)} \notag \\
& \qquad\quad \cdot \prod_{i=1}^\nu \qBin{n_i+iL-2\sum_{l=1}^i N_l-\alpha_{i,s}-\alpha_{i,b}}{n_i}{q}.
\label{eq:4.1}
\end{align}
Next, we sum over $\mathbf n$ to obtain
\begin{align}
& \sum_{\mathbf n}q^{N_1^2+\cdots+N_\nu^2+N_{s+1}+\cdots+N_\nu-M(N_1+\cdots+N_\nu)}
\prod_{i=1}^\nu \qBin{n_i+iL-2\sum_{l=1}^i N_l-\alpha_{i,s}-\alpha_{i,b}}{n_i}{q} \notag \\
& \qquad = \sum_{\mathbf n}\sum_{p\in\mathbb S_{-M,L-M}^\nu(\nu-s,\nu-b,\mathbf n)}q^{wg(p)}=
G_{-M,L-M}^\nu(\nu-s,\nu-b,q).
\label{eq:4.2}
\end{align}

\noindent
{\it Second step}. Here, the argument is exactly the same as the one used in
deriving \eqn{1.29}. We simply state the result
\begin{equation}
G_{-M,L-M}^\nu(\nu-s,\nu-b,q)=
\sum_{s'=0}^\nu G_{-M,0}^\nu(\nu-s,s',q)\; G_{0,L-M}^\nu(s',\nu-b,q).
\label{eq:4.3}
\end{equation}

\noindent
{\it Third step}. If the path $\tilde p\in\mathbb{\tilde S}_{-M,0}^\nu(s,b)$ is obtained 
from the path $p\in\mathbb{\tilde S}_{0,M}^\nu(b,s)$ by reflecting $p$ across the $y$-axis, 
then
\begin{equation}
wg(\tilde p)=-wg(p).
\label{eq:4.4}
\end{equation}
Hence,
\begin{align}
G_{-M,0}^\nu(\nu-s,s',q)& =
\sum_{p\in\tilde{\mathbb S}_{-M,0}^\nu(\nu-s,s')}q^{wg(p)} \notag \\
& = \sum_{p\in\tilde{\mathbb S}_{0,M}^\nu(s',\nu-s)}\Big(\frac{1}{q}\Big)^{wg(p)}=
G_{0,M}^\nu(s',\nu-s,\frac{1}{q}).
\label{eq:4.5}
\end{align}
Next, combining \eqn{4.2}, \eqn{4.3} and \eqn{4.5} we find that
\begin{align}
& \sum_{\mathbf n}q^{N_1^2+\cdots+N_\nu^2+N_{s+1}+\cdots+N_\nu-M(N_1+\cdots+N_\nu)}
 \prod_{i=1}^\nu \qBin{n_i+iL-2\sum_{l=1}^i N_l-\alpha_{i,s}-\alpha_{i,b}}{n_i}{q}\notag \\
& \qquad = \sum_{s'=0}^\nu G_{0,M}^\nu(s',\nu-s,\frac{1}{q})\; G_{0,L-M}^\nu(s',\nu-b,q).
\label{eq:4.6}
\end{align}
It follows from Theorem~\ref{thm:1} that
\begin{equation}
\lim_{L\rightarrow\infty}G_{0,L}^\nu(s,b,q)=
\frac{1}{\prod_{j\ge 1,j\not\equiv 0,\pm(\nu-s+1)(\mymod 2\nu+3)}(1-q^j)}.
\label{eq:4.7}
\end{equation}
Letting $L$ tend to infinity in \eqn{4.6}, we obtain with the aid of \eqn{1.24} and \eqn{4.7}
\begin{align}
& \sum_{\mathbf n}\frac{q^{N_1^2+\cdots+N_\nu^2+N_{s+1}+\cdots+N_\nu-M(N_1+\cdots+N_\nu)}}
{(q)_{n_1}(q)_{n_2}\ldots(q)_{n_\nu}} \notag \\
& \qquad = \sum_{s'=0}^\nu \frac{G_{0,M}^\nu(\nu-s',\nu-s,\frac{1}{q})}
{\prod_{j\ge 1,j\not\equiv 0,\pm(s'+1)(\mymod 2\nu+3)}(1-q^j)}.
\label{eq:4.8}
\end{align}
Finally, recalling \eqn{2.16} we arrive at the desired second variant
\begin{align}
& \sum_{\mathbf n}\frac{q^{N_1^2+\cdots+N_\nu^2+N_{s+1}+\cdots+N_\nu-M(N_1+\cdots+N_\nu)}}
{(q)_{n_1}(q)_{n_2}\ldots(q)_{n_\nu}} \notag \\
& \qquad = \sum_{\stackrel{s'=1}{s+s'\not\equiv\nu M(\mymod 2)}}^{2\nu+2}
\frac{W_{s',\nu-s}^\nu(M,\frac{1}{q})}
{\prod_{j\ge 1,j\not\equiv 0,\pm s'(\mymod 2\nu+3)}(1-q^j)}.
\label{eq:4.9}
\end{align}
Setting $M=1$ in \eqn{4.9} we easily obtain
\begin{equation}
\sum_{\mathbf n}\frac{q^{N_1^2+\cdots+N_\nu^2-N_1-N_2-\cdots-N_s}}
{(q)_{n_1}(q)_{n_2}\ldots(q)_{n_\nu}}=
\sum_{s'=\nu-s+1}^{\nu+1}
\prod_{\stackrel{j\ge 1}{j\not\equiv 0,\pm s'(\mymod 2\nu+3)}}\frac{1}{(1-q^j)},
\label{eq:4.10}
\end{equation}
which is essentially the same (modulo misprint) as identity~(3.3) in~\cite{B2}.

\section{Further generalizations}
\label{sec:5}

\noindent
Comparing \eqn{1.32} and \eqn{4.9} with $s=\nu$, the structural similarities of these two 
formulas become obvious. This resemblance strongly suggests that our two
variants are special cases of a more general formula for multisums of the form
\begin{equation}
\sum_{\mathbf n}\frac{q^{N_1^2+\cdots+N_\nu^2-M_1N_1-\cdots-M_\nu N_\nu}}
{(q)_{n_1}(q)_{n_2}\ldots(q)_{n_\nu}},
\label{eq:5.1}
\end{equation}
with $(M_1,M_2,\ldots,M_\nu)\in\mathbb Z^\nu$.

\noindent
The unifying formula requires generalized $q$-multinomials that depend on $\nu$
finitization parameters. Fortunately, such objects have already appeared in the literature
\cite{BU,SW}:
\begin{align}
\qBin{\mathbf L}{a}{\nu}:=
\sum_{j_1+\cdots+j_\nu=a+\frac{1}{2}\sum_{i=1}^\nu iL_i} &
q^{\sum_{l=2}^\nu j_{l-1}(L_l+\cdots+L_\nu-j_l)} \notag\\
& \cdot \qBin{L_\nu}{j_\nu}{q}\qBin{L_{\nu-1}+j_\nu}{j_{\nu-1}}{q}\cdots\qBin{L_1+j_2}{j_1}{q},
\label{eq:5.2}
\end{align}
where
\begin{equation}
\mathbf L=(L_1,L_2,\ldots,L_\nu)\in\mathbb Z^\nu,
\label{eq:5.3}
\end{equation}
and
\begin{equation}
a+\frac{1}{2}\sum_{i=1}^\nu iL_i\in\mathbb Z_{\ge 0}.
\label{eq:5.4}
\end{equation}
In \cite{SW}, the polynomials (5.2) were termed $q$-supernomial coefficients.

\noindent
We notice that
\begin{equation}
\qBin{(L,0,\ldots,0)}{a}{\nu}=\qBin{L}{a+\frac{L}{2}}{q},
\label{eq:5.5}
\end{equation}
and
\begin{equation}
\qBin{(0,\ldots,0,L)}{a}{\nu}=\qBinpwr{L}{a+\frac{\nu L}{2}}{\nu}{0}.
\label{eq:5.6}
\end{equation}
Using \eqn{5.5}, \eqn{5.6} along with \eqn{1.32} and \eqn{4.9} with
$s=\nu$, it is easy to guess that
\begin{align}
& \sum_{\mathbf n}\frac{q^{N_1^2+\cdots+N_\nu^2-M_1N_1-\cdots-M_\nu N_\nu}}
{(q)_{n_1}(q)_{n_2}\ldots(q)_{n_\nu}} \notag\\
& \qquad =\sum_{\stackrel{s=1}{s+\nu+\sum_{i=1}^\nu M_i\mbox{\scriptsize{ odd}}}}^{2\nu+2}
\frac{I_{s,\nu+1}^\nu (\tilde{\mathbf{M}},\frac{1}{q})}
{\prod_{j\ge 1,j\not\equiv 0,\pm s(\mymod 2\nu+3)}(1-q^j)},
\label{eq:5.7}
\end{align}
where
\begin{equation}
\tilde{\mathbf M}=(M_1-M_2,M_2-M_3,\dots,M_{\nu-1}-M_\nu,M_\nu),
\label{eq:5.8}
\end{equation}
and
\begin{align}
I_{s,b}^\nu(\mathbf L,q) & := \sum_{j=-\infty}^\infty \Bigg\{
q^{j((2j+1)(2\nu+3)-2s)}
\qBin{\mathbf L}{\frac{b-s}{2}+(2\nu+3)j}{\nu} \notag \\
& \qquad\qquad - q^{(2j+1)((2\nu+3)j+s)}
\qBin{\mathbf L}{\frac{b+s}{2}+(2\nu+3)j}{\nu} \Bigg\}.
\label{eq:5.9}
\end{align}
The remarkable formula \eqn{5.7} and its finite analogs will be the subject of our next paper. 

\subsection*{Acknowledgment}

We would like to thank G.E.~Andrews and C.~Krattenthaler for their interest.
We are grateful to A.~Riese for his help in preparing the figures for this manuscript, and
S.O.~Warnaar for many patient explanations of his work \cite{W} and for his comments.


\end{document}